\documentclass{article}

\usepackage{graphicx}
\usepackage{amsmath}
\usepackage{amssymb}
\usepackage[mathscr]{eucal}

\def \calB{{\cal B}}

\def \calE{{\cal E}}
\def \calF{{\cal F}}

\def \PointP{{\mathscr{P}}}

\def \mone{{\text{-}1}}

\def\mumu{{\raisebox{1 ex}{$_{_\mu}$}}}
\def\mumumu{{\hskip -1.95 mm\mumu}}
\def\nunu{{\raisebox{1 ex}{$_{_\nu}$}}}
\def\nununu{{\hskip -1.84 mm\nunu}}

\def\bibref{\par\noindent\hangindent=20pt}

\title{Image and Reciprocal Image of a Measure. Compatibility Theorem.}
\author{Albert Tarantola\thanks{Institut de Physique du Globe de Paris, 4 place Jussieu, 75005 Paris, France. E-mail: tarantola@ipgp.jussieu.fr}}

\begin{document}
\maketitle

%
%

\bigskip\bigskip\bigskip\bigskip

\begin{abstract}
It is proposed that to the usual probability theory, three definitions and a new theorem are added, the resulting theory allows one to displace the central role usually given to the notion of conditional probability. When a mapping \,$ \varphi $\, is defined between two measurable spaces, to each measure \,$\mu$\, introduced on the first space, there corresponds an image \,$ \varphi[\mu] $\, on the second space, and, reciprocally, to each measure \,$ \nu $\, defined on the second space the corresponds a reciprocal image \,$ \varphi^\mone [\nu] $\, on the first space. As the intersection \,$ \cap $\, of two measures is easy to introduce, a relation like \,$ \varphi[\, \mu \,\cap\, \varphi^\mone [\nu] \,]= \varphi[\mu] \,\cap\, \nu $\, makes sense. It is, indeed, a theorem of the theory. This theorem gives mathematical consistency to inferences drawn from physical measurements.
\end{abstract}

%
%

\newpage
\quad
\bigskip\bigskip\bigskip
\tableofcontents
\newpage

\newpage

%
%

\section{Preliminary}

Assume given a measurable space%
\footnote{%
As usual, here \,$ \Omega $\, denotes a set and \,$ \calF $\, is a collection of subsets of \,$ \Omega $\, that is a $\sigma$-field (i.e., \,$ \calF $\, is nonempty and it is closed under complementation and countable unions of its members).}%
%
%
\,$ (\Omega,\calF)$\,, and a measure%
\footnote{%
A (positive) measure (measures are implicitly assumed to be positive) is a function \,$ \mu : \calF \mapsto [0,\infty] $\, satisfying two properties: 
(\emph{i}) the measure of the empty set is zero, and  
(\emph{ii}) the measure of the union of a countable sequence of pairwise disjoint sets in \,$ \calF $\, equals the sum of the measures of each set.}%
%
%
\,$ \mu $\, on 
\,$ (\Omega,\calF) $\, that is $\sigma$-finite%
\footnote{%
A (positive) measure \,$ \mu $\, defined on a $\sigma$-algebra \,$ \calF $\, of subsets of a set \,$ \Omega $\, is called \emph{$\sigma$-finite} if \,$ \Omega $\, is the countable union of measurable sets of finite $\mu$-measure. A~set in a measure space has \emph{$\sigma$-finite measure} if it is a countable union of sets with finite measure.
A measure \,$ \mu $\, is called \emph{finite} if \,$ \mu[\Omega] $\, is a finite real number (rather than $ \infty $). A $\sigma$-finite measure may not be finite (the Lebesgue measure on the real line is $\sigma$-finite, but not finite). An example of a measure on the real line that is not $\sigma$-finite is the \emph{counting measure} (the counting measure of a set of real numbers is the number of elements in the set): every set with finite measure contains only finitely many real numbers, and it would take uncountably many such sets to cover the entire real line. The Radon-Nikodym theorem does not apply to the counting measure and  no density can be associated to it. It is to avoid such ``pathological'' measures that the $\sigma$-finite hypothesis is introduced.}.
%
%
Then, \,$ (\Omega,\calF,\mu)$\, is called a $\sigma$-finite measure space. Let \,$ \nu $\, be a second measure on \,$ (\Omega,\calF)$\,. The following assertions ---a symmetric version of the Radon-Nikodym theorem--- are equivalent (Schilling, 2006). 
\begin{itemize}

\item
The measure \,$ \nu $\, is absolutely continuous%
\footnote{
The measure \,$ \nu $\, is said to be \emph{absolutely continuous} with respect to the measure \,$ \mu $\, if
\,$ \mu[F] = 0 \ \Rightarrow \ \nu[F] = 0 $\,. One writes \,$ \nu \ll \mu $\,.}
%
%
with respect to \,$ \mu $\,.

\item
There is a $\mu$-almost everywhere unique function from \,$ \Omega $\, into \,$ [0,\infty) $\,, denoted \,$ d\nu/d\mu $\,, such that
\begin{equation}
\nu[F] \ = \ \int_F \frac{d\nu}{d\mu} \ d\mu \qquad \text{for every} \ F \in\calF \quad . 
\end{equation}
\end{itemize}
The function \,$ d\nu/d\mu $\, is called the Radon-Nikodym \emph{density} associated with \,$ \nu $\, by~\,$ \mu $\,, or the Radon-Nikodym \emph{derivative} of \,$ \nu $\, with respect to \,$ \mu $\, .

If a $\sigma$-finite measure \,$ \mu $\, is such that \,$ \mu[\Omega] = 1 $\,, one says that \,$ \mu $\, is a \emph{probability measure}, and the measure \,$ \mu[F] $\, of some set \,$ F \in \calF $\, is then called the probability of the set%
\footnote{
When dealing with probability measures only, sets are often called \emph{events}.}
%
%
\,$ F $\,. 

%
%

\section{Definitions and properties}

%
%

\subsection{Intersection of measures}

Given a $\sigma$-finite measure space \,($\Omega,\calF,\mu$)\,, consider two $\sigma$-finite measures \,$ \nu_1 $\, and \,$ \nu_2 $\,, at least one of them ---\,say \,$ \nu_1 $\,--- being absolutely continuous with respect to the base measure \,$ \mu $\,.

\emph{
{\bf Definition:}
Given some finite constant \,$ n $\,, the \emph{intersection} of the two measures \,$ \nu_1 $\, and \,$ \nu_2 $\,, is the measure denoted \,$ \nu_1 \cap \nu_2 $\, and defined as
\begin{equation}
\boxed{\qquad
(\nu_1 \cap \nu_2)[F] \ = \ \frac{1}{n} \int_F \frac{d\nu_1}{d\mu} \ d\nu_2 \qquad 
\text{for every} \ \ F \in \calF \quad . \quad }
\label{eq: intersec_3822}
\end{equation}
}%
\noindent
It is obvious that this defines a measure and that ---by virtue of the Radon-Nikodym theorem--- it is absolutely continuous w.r.t.\ \,$ \nu_2 $\,.
The operation \,$ \cap $\, depends on the base measure \,$ \mu $\,, so, when necessary, a more explicit notation, like \,$ \cap_\mumumu $\,, can be used.

Remark: When dealing with arbitrary measures, one may well take \,$ n = 1 $\,, while when dealing with probability measures, it may be more convenient to take
\,$ n \, = \, \int_\Omega (d\nu_1/d\mu) \ d\nu_2 $\,
as, then, \,$ (\nu_1\cap\nu_2)[\Omega] = 1 $\,, this implying that the intersection of two probability measures is a probability measure (but, then, the intersection would only be defined if \,$ 0 < n < \infty $\,).

Should the measure \,$ \nu_2 $\, also be absolutely continuous with respect to \,$ \mu $\,, 
equation~\eqref{eq: intersec_3822} could be written
\,$ (\nu_1 \cap \nu_2)[F] \, = \, \frac{1}{n} \int_F \frac{d\nu_1}{d\mu} \, \frac{d\nu_2}{d\mu} \, d\mu $\,,
the measure \,$ \nu_1 \cap \nu_2 $\, would also be absolutely continuous with respect to \,$ \mu $\,, and its Radon-Nikodym density would be
\begin{equation}
\frac{d(\nu_1\cap\nu_2)}{d\mu} \ = \ \frac{1}{n} \, \frac{d\nu_1}{d\mu} \, \frac{d\nu_2}{d\mu} \quad .
\label{eq: fautyaller_92773}
\end{equation}

Comment: The term {\sl intersection} is justified in section~\ref{sec: Measures versus sets: intersection}, when the intersection of (measurable) sets is found as a special instance of the intersection of measures. Another special instance of the intersection of (probability) measures corresponds to the notion of conditional probability (see section~\ref{Intersection of measures and conditional probability}). 

%
%

\subsection{Reciprocal image of a measure}

Let \,$ (X,\calE) $\, and \,$ (Y,\calF) $\, be two measurable spaces, and \,$ \varphi : X \mapsto Y $\, a measurable%
\footnote{
The mapping \,$ \varphi $\, is measurable, if the reciprocal image of every set in \,$ \calF $\, is in \,$ \calE $\,. Non measurable mappings are generally considered pathological.}
%
%
mapping. 
Two measures \,$ \mu $\, and \,$ \nu $\, are introduced (to be considered as base measures) such that \,$ (X,\calE,\mu) $\, and \,$ (Y,\calF,\nu) $\, are $\sigma$-finite measure spaces.

\emph{
{\bf Definition:}
Given some finite constant \,$ n $\,, to every measure \,$ \tau $\, on \,$ (Y,\calF) $\, that is absolutely continuous with respect to \,$ \nu $\,, is associated a measure on \,$ (X,\calE) $\,, called the \emph{reciprocal image} of \,$ \tau $\,, denoted \,$ \varphi^\mone[\,\tau\,] $\,, and defined via
\begin{equation}
\boxed{ \qquad 
\frac{d(\varphi^\mone[\,\tau\,])}{d\mu} \ = \ \frac{1}{n} \ \Big( \frac{d\tau}{d\nu} \circ \varphi \Big) \quad . \quad }
\label{eq: erecp_46}
\end{equation}
}

Then, for every \,$ E \in \varphi^\mone[\calF] \subseteq \calE $\,, one has%
\footnote{
As \,$ \varphi $\, is a measurable mapping, and the function \,$ d\tau/d\nu $\, is measurable (with respect to \,$ \calF $\,), the function \,$ (d\tau/d\nu)\circ \varphi $\, is measurable (with respect to \,$ \varphi^\mone[\calF] \subseteq \calE $\,). (See, e.g., Halmos, 1950).}
%
%
\,$ (\varphi^\mone[\,\tau\,])[E] \, = \, \frac{1}{n} \int_E \frac{d\tau}{d\nu}(\varphi(x)) $ $ d\nu(x) $\,, and the Radon-Nikodym theorem ensures that \,$ \varphi^\mone[\,\tau\,] $\, is, indeed, a measure. As this reciprocal image depends on the two base measures, the more explicit notation \,$ \varphi^\mone[\,\tau\,;\,\mu,\nu\,] $\, can be used.

Remark: When dealing with arbitrary measures, one may well take \,$ n = 1 $\,,
while when dealing with probability measures, it may be more convenient to take
\,$ n = \int_X \frac{d\tau}{d\nu} \circ \varphi \ d\mu $\,,
as, then, \,$  (\varphi^\mone[\tau])[X] = 1 $\,, this implying that the reciprocal image of a probability measure is a probability measure (but, then, the reciprocal image would only be defined if \,$ 0 < n < \infty $\,).

%
%

\subsection{Image of a measure}
\label{sec: Image of a measure}

Let \,$ (X,\calE )$\, and \,$ (Y,\calF) $\, be two measurable spaces, and \,$ \varphi : X \mapsto Y $\, a measurable mapping. 

\emph{
{\bf Definition:}
To every measure \,$ \pi $\, on \,$ (X,\calE) $\,, is associated a measure%
\footnote{
It is not difficult to verify that \,$ \varphi[\pi] $\, is, indeed, a measure.
First, the measure of the empty set is zero, because, by definition of the reciprocal
mapping,\,$ \varphi^\mone [\,\emptyset\,] = \,\emptyset\, $\,, so \,$
(\varphi[\pi])[\,\emptyset\,] = \pi[\, \varphi^\mone [\,\emptyset\,] \,] =
\pi[\,\emptyset\,] = 0 $\,. Second, we have to check that if \,$ F_1,
F_2,\dots $\, is a countable sequence of pairwise disjoint sets in \,$
\calF $\, the measure of the union of all the \,$ F_i $\, is equal to
the sum of the measures of each \,$ F_i $\,: \,$ (\varphi[\pi])[\,
\bigcup_i F_i \, ] \, = \, \sum_i \tau[F_i] $\,. First, by definition 
of image of a measure,
one has \,$ (\varphi[\pi])[\, \bigcup_i F_i \, ]  =  \pi[ \, \varphi^\mone [\,
\bigcup_i F_i \, ] \, ] $\, and, as the reciprocal image of a union is
the union of the reciprocal images, \,$ (\varphi[\pi])[\, \bigcup_i F_i \, ] \, = \,
\pi[ \, \bigcup_i \varphi^\mone [F_i] \, ] $\,. But \,$ \pi $\, is a
measure, and the reciprocal image of disjoint sets is disjoint, so
\,$(\varphi[\pi])[\, \bigcup_i F_i \, ] \, = \, \sum_i \pi[ \, \varphi^\mone [F_i]
\, ] $\,. Finally, using again the definition of image of a measure,
this leads to desired property. We have thus checked that the image of a 
measure is a measure.}
%
%
on \,$ (Y,\calF) $\,, denoted \,$ \varphi[\pi] $\,, called the \emph{image measure}:
\begin{equation}
\boxed{ \qquad 
\varphi[\pi] \ = \ \pi \circ \varphi^\mone \quad , \quad }
\label{eq: notveryfarfrom_993}
\end{equation}
i.e., explicitly,
\,$ (\varphi[\pi])[F] \, = \, \pi[ \, \varphi^\mone[F] \, ] $\, 
for every \,$ F \in \calF $\,.
}

The measure \,$ \varphi[\pi] $\, needs not be%
\footnote{
As an example, this happens when \,$ X = \Re^p $\, and \,$ Y = \Re^q $\, with \,$ p < q $\,, and a continuous mapping \,$ \varphi $\, (with the standard Lebesgue measures assumed), because then  \,$ \varphi[X] $\, is a $p$-dimensional submanifold of \,$ \Re^q $\,.}
%
%
absolutely continuous with respect to some base measure, so \,$ \varphi[\pi] $\, may not be representable by a bona-fide density%
\footnote{
When \,$ \varphi[\pi] $\, is representable by a density, it is, in general,
easy to find an expression of it, but the (elementary) methods
to be used are quite different in every situation
(see examples in sections~\ref{sec: Mapping between discrete sets} and~\ref{sec: Propagation of uncertainties in physical measurements}), 
and a general expression for the density is not available.}.
%
%
This does not cause any complication in the applications we have in mind.

We shall later need the following property (Halmos, 1950, page 163): for any measurable function \,$ K $\, and any set \,$ F \in \calF $\,,
\begin{equation}
\int_F K \ d(\pi\circ\varphi^\mone) \ = \ \int_{\varphi^\mone[F]} K\circ \varphi \ \ d\pi \quad ,
\label{eq: to-dem-3339277}
\end{equation}
i.e., \,$ \int_F K(y) \ d(\varphi[\pi])(y) \ = \ \int_{\varphi^\mone[F]} K(\varphi(x)) \ d\pi(x) $\,.

\emph{Comment:} To have an intuitive idea of the notion ``image of a measure'', consider a collection of elements \,$ x_1, x_2, x_3 ,\dots $\, of \,$ X $\, that are independent sample elements of the measure \,$ \pi $\,. Then, it is easy to see that the elements \,$ \varphi(x_1), \varphi(x_2), \varphi(x_3) ,\dots $\, of~\,$ Y $\, are independent sample elements of the measure \,$ \varphi[\pi] $\,. In fact, this property alone may suggest introducing the notion of an image of a measure.

%
%

\subsection{Compatibility property}

Let \,$ (X,\calE,\mu) $\, and \,$ (Y,\calF,\nu) $\, be two $\sigma$-finite measure spaces, and \,$ \varphi : X \mapsto Y $\, be a measurable mapping.
Let \,$ \pi $\, be a measure over \,$ (X,\calE) $\, that is $\sigma$-finite, and \,$ \tau $\, a measure over  \,$ (Y,\calF) $\,, that is absolutely continuous with respect to the base measure~\,$ \nu $\,. 

\medskip

\emph{
{\bf Theorem:}
One always has
\begin{equation}
\boxed{\qquad
\varphi[ \, \pi \cap \pi' \, ] \ = \ \tau' \cap \tau \qquad 
\text{where} \qquad 
\begin{cases}
\quad \pi' \ = \ \varphi^\mone[\,\tau\,] \\
\quad \tau' \ = \ \varphi[\pi] \quad . \\
\end{cases}
\quad}
\label{eq: forisaad}
\end{equation}
}%

Note that while the measure \,$ \tau $\, is assumed to be absolutely continuous with respect to the base measure \,$ \nu $\,, the measures \,$ \varphi[\pi] $\, and \,$ \varphi[\pi\cap\pi'] $\, may be singular.

\medskip

To demonstrate the identity in equation~\eqref{eq: forisaad} means to verify that for any set \,$ F \in \calF $\,, one has
\,$ ( \, \varphi[ \, \pi \cap \pi' \, ] \, )[F] \, = \, (\,\tau' \cap \tau\,) [F] $\,.
This is done by writing the following sequence of identities (that successively use 
equations~\eqref{eq: notveryfarfrom_993},
\eqref{eq: intersec_3822},
\eqref{eq: erecp_46},
\eqref{eq: to-dem-3339277},
and~\eqref{eq: intersec_3822} again):
\begin{equation}
\begin{split}
( \, \varphi[ \, \pi \cap_\mumumu \varphi^\mone[\tau;\mu,\nu] \, ] \, )[F] \ 
& = \ (\, \pi \cap_\mumumu \varphi^\mone[\tau;\mu,\nu] \,)[\,\varphi^\mone[F]\,] \\
& = \ \frac{1}{n} \int_{\varphi^\mone[F]} \frac{d(\varphi^\mone[\tau;\mu,\nu])}{d\mu}(x) \ d\pi(x) \\
& = \ \frac{1}{n'} \int_{\varphi^\mone[F]} \frac{d\tau}{d\nu}(\varphi(x)) \ d\pi(x) \\
& = \ \frac{1}{n'} \int_F \frac{d\tau}{d\nu} (y) \ d(\varphi[\pi])(y) \\[2 pt]
& = \ \frac{1}{n''} \, ( \, \varphi[\pi] \cap_\nununu \tau \, )[F] \quad , \\
\end{split}
\end{equation}
so the property holds%
\footnote{
The constant \,$ n'' $\, equals one, because for general measures the two constants in 
equations~\eqref{eq: intersec_3822} and~\eqref{eq: erecp_46} should be taken equal to one,
while for probability measures, there is an automatic renormalization.}.
%
%

%
%

\section{Measures and sets}
\label{sec: Measures and sets}

%
%

\subsection{Measure-sets}

The definitions and properties above have a direct relation with definitions and properties in set theory, and, in some sense, they generalize them. To see this, let us start by introducing the notion of measure-set.

Let \,$ (\Omega,\calF,\mu) $\, be a $\sigma$-finite measure space.
To every set \,$ A \in \calF $\, we shall associate a measure, denoted \,$ \mu_A $\,, and defined via the condition
\begin{equation}
\mu_A[F] \ = \ \frac{1}{n_A} \ \mu[A\cap F] \qquad \text{for every} \ \ F \in \calF \quad ,
\end{equation}
where \,$ n_A $\, is a suitable chosen constant (that may depend on \,$ A $\, but not on~\,$ F $\,). 
As suggested above, one may well take
\begin{equation}
n_A \ = \ 1 \qquad \text{for arbitrary measures} \quad ,
\end{equation}
or
\begin{equation}
n_A \ = \ \mu[A] \qquad \text{for probability measures} \quad ,
\end{equation}
because, then, \,$ \mu_A[\Omega] = 1 $\, (of course, this assumes \,$ \mu[A] \neq 0 $\,.)
Such a measure shall be called a \emph{measure-set}, so we can talk about the measure-set \,$ \mu_A $\, associated with a set \,$ A $\, by a measure \,$ \mu $\,.
The $\mu$-density associated with a measure-set \,$ \mu_A $\, is clearly%
\footnote{
As, for every \,$ F \in \calF $\,, \,$ \frac{1}{n_A} \int_F \chi_A \, d\mu \, = \, 
\frac{1}{n_A} \int_{A\cap F} d\mu \, = \, \frac{1}{n_A} \, \mu(A \cap F) \, = \, \mu_A[F] $\,.}
%
%
\,$ \frac{1}{n_A} \, \chi_A $\,, i.e., proportional to the characteristic function of the set~\,$ A $\,.

Of course, there may be subsets of \,$ \Omega $\, that are not in \,$ \calF $\,, but, as far as one is only interested in the sets in \,$ \calF $\,, one can consider that any measure absolutely continuous 
w.r.t.\ \,$\mu $\, is something like a generalized set: while a measure-set can be identified to a set (its density taking only two possible values, \,$ 0$\, or \,$1/n_A $\,), an arbitrary measure (with a density taking any nonnegative value) is a kind of generalized object, that contains measure-sets and, therefore, sets as special cases.

The names given to the three notions introduced above ---intersection of measures, and image and reciprocal image of a measure--- are justified because (as we are about to see) when applied to measure-sets they do correspond to the intersection, the image and the reciprocal image of sets.

%
%

\subsection{Measures versus sets: intersection}
\label{sec: Measures versus sets: intersection}

If \,$ A $\, and \,$ B $\, are two sets in~\,$ \calF $\, and \,$ \mu_A $\, and \,$ \mu_B $\, are the two associated measure-sets, one has%
\footnote{%
For any \,$ F \in \calF $\,, 
\,$ (\mu_A \cap \mu_B)[F] = \frac{1}{n_A n_B} \int_F \chi_A \, \chi_B \, d\mu 
= \frac{1}{n_A n_B} \int_{F\cap A\cap B} d\mu 
= \frac{n_{A\cap B}}{n_A n_B} \int_F \chi_{A\cap B} d\mu 
= \frac{n_{A\cap B}}{n_A n_B} \ \mu_{A\cap B}[F]
$\,. }
%
%
\begin{equation}
\mu_A \cap \mu_B \ = \ k \, \mu_{A \cap B} 
\end{equation}
with the constant%
\footnote{
For general measures, \,$ k = 1 $\,, while for probability measures, \,$ k = \frac{\mu[A\cap B]}{\mu[A] \mu[B]} $\,.}
%
%
\,$ k = \frac{n_{A\cap B}}{n_A \, n_B} $\,.
So ---in the special case where the measures are measure-sets--- the definition of intersection of measures is consistent with the definition of intersection of sets.

%
%

\subsection{Intersection of measures and conditional probability}
\label{Intersection of measures and conditional probability}

Letting \,$ A $\, be a fixed set of \,$ \calF $\,, let us now consider the intersection of an arbitrary measure \,$ \nu $\, and the measure-set \,$ \mu_A $\,, i.e., the measure \,$ \nu \cap \mu_A $\,.
One has%
\footnote{
For
\,$ (\nu \cap \mu_A)[F] \propto \int_F \chi_A \, d\nu
= \int_{F\cap A} d\nu = \nu[F\cap A] $\,}
%
%
\begin{equation}
(\nu \cap \mu_A)[F] \ = \ \frac{1}{n} \ \nu[F\cap A] \qquad \text{for every} \ \ F \in \calF \quad ,
\end{equation}
where \,$ n $\, is a constant.
This is particularly interesting when dealing with probability measures, because, then, \,$ n = \nu[A] $\,,
and one has
\begin{equation}
(\nu \cap \mu_A)[F] \ = \ \frac{\nu[F\cap A]}{\nu[A]} \qquad \text{for every} \ \ F \in \calF \quad .
\end{equation}
One immediately recognizes there the expression of Kolmogorov's conditional probability, usually denoted
\,$ \nu[F\cap A] / \nu[A] = \nu[ \, F \, | \, A \, ] $\,. Using this notation,
\begin{equation}
(\nu \cap \mu_A)[F] \ = \ \nu[ \, F \, | \, A \, ] \qquad \text{for every} \ \ F \in \calF \quad .
\end{equation}
So we have the following

\emph{
{\bf Property:}
For every given probability measure \,$ \nu $\,, Kolmogorov's conditional probability, given any set \,$ A $\,, is identical to the intersection of \,$ \nu $\, by the measure-set \,$ \mu_A $\,.}

So we see that the notion of conditional probability is a special case of the notion of intersection of measures: when evaluating the intersection of an arbitrary measure by a measure set, we have the conditional probability, but we can evaluate the intersection of two general measures. I claim that there are problems that are naturally formulated in terms of the intersection of two measures (see 
sections~\ref{sec: Intersection of probability measures}
and~\ref{sec: Interpretation of observations (2: using measures)}
for examples). As this notion has not been available so far, some hand-waving has been necessary to make this kind of problems fit into the available mathematical structure. This, plus the fact that general mappings between arbitrary sets (as opposed to linear mappings between linear spaces) can be taken as root elements, is what has motivated the building of the present theory.

%
%

\subsection{Measures versus sets: reciprocal images}

When considering a mapping \,$ \varphi $\, from a set \,$ X $\, into a set \,$ Y $\, to every set \,$ B \subseteq Y $\, there is associated a subset of \,$ X $\,, denoted \,$ \varphi^\mone [B] $\, and named the reciprocal image%
\footnote{
The set \,$ \varphi^\mone [B] $\, is made of all the elements of \,$ X $\, whose image is in \,$ B $\,.}
%
%
of the set~\,$ B $\,. But the reciprocal image of a set can also be defined in terms of the characteristic functions of the sets: letting \,$ \xi_A $\, the characteristic function of a set \,$ A \subseteq X $\, and
\,$ \chi_B $\, that of a set \,$ B \subseteq Y $\,, one clearly has
\begin{equation}
\xi_{\varphi^\mone[B]} \ = \ \chi_B \circ \varphi \quad ,
\end{equation}
a relation that is identical (with \,$ n = 1 $\,) to the relation~\eqref{eq: erecp_46} expressing the reciprocal image of a measure in terms of Radon-Nikodym densities.
So, as it already happened with the intersection of measures, the notion of reciprocal image of a measure is consistent with the definition of reciprocal image of a set: the set associated to the reciprocal image of the measure-set that is associated to a set \,$ B $\, is the reciprocal image of the set \,$ B $\,:
\begin{equation}
\varphi^\mone [\nu_B] \ = \ \mu_{\varphi^\mone [B]} \quad . 
\end{equation}
In this sense, again, the notion of reciprocal image of a set is ``contained'' in the notion of reciprocal image of a measure.

%
%

\subsection{Measures versus sets: images}

The relation between the notion of image of a set (in set theory) and the notion of image of a measure is subtle. In this short note, let us just mention that the support%
\footnote{
The support of a function is the set of points where the function is not zero.}
%
%
of the image of a measure \,$ d(\varphi[\pi])/d\nu $\, is the image (in the sense of set theory) of the support of the original  measure \,$ d\pi/d\mu $\,.

%
%

\subsection{Measures versus sets: compatibility property}

In set theory, for arbitrary sets \,$ A_1 $\, and \,$ A_2 $\, and for an arbitrary mapping \,$ \varphi $\,, one has
\begin{equation}
\varphi[ \, A_1 \cap A_2  \,] \subseteq \varphi[A_1] \cap \varphi[A_2] \quad ,
\end{equation}
a relation that is well-known but not very useful here. A more useful relation (for making inferences involving sets and mappings) is that, for arbitrary sets \,$ A $\, and \,$ B $\,, and an arbitrary mapping \,$\varphi $\,, one has
\begin{equation}
\varphi[ \, A \cap \varphi^\mone [B] \, ] \ = \ \varphi[A] \cap B \quad .
\label{eq: forisaad_332}
\end{equation}
For reasons that shall become clear in the applications (see section~\ref{sec: Interpretation of observations (2: using measures)}), it is interesting to extend this identity into probability theory (or, more generally, inside measure theory). But, of course, our compatibility property (equation~\eqref{eq: forisaad})
\begin{equation}
\varphi[ \, \pi \cap \varphi^\mone[\,\tau\,] \, ] \ = \ \varphi[\pi] \cap \tau \quad ,
\end{equation}
is identical to relation~\eqref{eq: forisaad_332}, excepted that it concerns measures instead of sets. So, in some sense,
we have generalized the set relation. In any case, when the relation \,$ \varphi[ \, \pi \cap \varphi^\mone[\,\tau\,] \, ] \, = \, \varphi[\pi] \cap \tau $\, is applied to measure sets, it becomes relation~\eqref{eq: forisaad_332}.

%
%

\section{Applications}
\label{sec: Applications}

%
%
\subsection{Intersection of probability measures}
\label{sec: Intersection of probability measures}

Let \,$ \Omega $\, represent the surface of the sphere 
of unit radius, and \,$ \calB $\, the usual Borel collection of subsets%
\footnote{%
The Borel sigma-field is defined as the smallest sigma-field containing all the open subsets.}.
%
%
Consider, on the measurable space \,$ (S,\calB) $\,, the ordinary surface 
measure: for any set \,$ B \in \calB$\, of points on the sphere,
\,$ S[B] $\, is the surface of \,$ B $\,. Two probability measures 
\,$ P_1 $\, and \,$ P_2 $\, are then considered, and
two points \,$ \PointP_1 $\, and \,$ \PointP_2 $\, are randomly created%
\footnote{%
The notion of random point is not introduced here; it is assumed that reader knows the basic notion of sampling from a probability measure.}
%
%
on the surface of the sphere, that are random point samples of the respective probability measures \,$ P_1 $\, and \,$ P_2 $\,.
If \,$ \PointP_1 \neq \PointP_2 $\, the two points are discarded, and two new points are generated. And so on until the two points happen to be identical, \,$ \PointP_1 = \PointP_2 = \PointP $\,\,.

\medskip

\emph{Question:} of which probability measure is \,$ \PointP $\,\, a random point sample?

\medskip

\emph{Answer:} point \,$ \PointP $\, is a random point sample of the probability measure
\begin{equation}
P \ = \ P_1 \cap P_2 \quad .
\label{eq: intersph_49}
\end{equation}


\emph{Proof:} The probability that the two points \,$ \PointP_1 $\, and \,$ \PointP_2 $\, happen to be identical is zero, so the question makes no immediate sense, and needs to be slightly reformulated. If the sphere is assumed to be tiled with a finite collection of (spherical) tiles of identical surface \,$ \Delta S $\,, then it can happen that the two points \,$ \PointP_1 $\, and \,$ \PointP_2 $\, are in the same tile. The finite probability that this happens in a given tile can then be evaluated (it is the renormalized product of the two probabilities assigned to the tile by each of the two probability measures \,$ P_1 $\, and \,$ P_1 $\,), and it is when taking the limit \,$ \Delta S \to 0 $\, that one gets the result. 

\medskip

Introducing the three probability densities \,$ f_1 $\,, \,$ f_2 $\,, and \,$ f_1\cap f_2 $\, associated with the three probability measures \,$ P_1 $\,, \,$ P_2 $\,, and \,$ P_1\cap P_2 $\, via
\begin{equation}
P_1[B] \, = \int_B f_1 \, dS \ \ ; \ \
P_2[B] \, = \int_B f_2 \, dS \ \ ; \ \
(P_1\cap P_2)[B] \, = \int_B (f_1\cap f_2) \, dS \ \ ,
\end{equation}
gives here, using equation~\eqref{eq: fautyaller_92773},
\begin{equation}
(f_1\cap f_2) \ = \ \frac{f_1 \, f_2}{\int_\Omega f_1 \, f_2 \, dS} \quad .
\label{eq: almodo_00922}
\end{equation}

To pass from this purely mathematical exercise to a problem involving real-life measurements, assume that two totally ``disentangled''%
\footnote{
I am trying here to avoid the use of the term \emph{independent} that has a related ---but different--- connotation in probability theory.}
%
%
measurements of the position of a floating object on the ocean provide the information described (following ISO's recommendations [ISO, 1993]) by two probability densities \,$ f_1 $\, and \,$ f_2 $\,. How should they be ``combined'' to represent the total available information? The detailed justifications of this is outside the scope of this short note, but I suggest here that experimental uncertainties are defined in such a way (ISO's way) that the answer to the question is precisely that in equation~\eqref{eq: almodo_00922}.

%
%

\subsection{Mapping between discrete sets}
\label{sec: Mapping between discrete sets}

Let \,$ X $\, and \,$ Y $\, be discrete spaces, \,$ \calE $\, and \,$ \calF $\, the respective collections of their subsets, and \,$ \varphi $\, a mapping from \,$ X $\, into \,$ Y $\,.
Consider that a probability measure \,$ \pi $\, on \,$ (X,\calE) $\,, is sampled, this providing elements \,$ x_1, x_2, \dots  $\, of \,$ X $\,, and therefore, via the mapping \,$ \varphi $\,, the image elements \,$ y_1 = \varphi(x_1) \, , \, y_2 = \varphi(x_2) \, , \, \dots $\, of \,$ Y $\,. Of which probability measure \,$ \tau $\, on \,$ (Y,\calF) $\, are the elements \,$ y_1,y_2,\dots $\, sample points?

The answer is \,$ \tau = \varphi[\pi] $\,, as this clearly corresponds to the very definition of image of a measure (equation~\eqref{eq: notveryfarfrom_993}):
\begin{equation}
\tau[F] \ = \ \pi[\, \varphi^\mone [F] \,] \qquad \text{for every} \ \ F \in \calF \quad .
\label{eq: image_after_halmos_3962}
\end{equation}

To transform this result into an explicit expression, we can introduce two base measures \,$ \mu $\, and \,$ \nu $\, on \,$ (X,\calE) $\, and \,$ (Y,\calF) $\, respectively, for instance, the respective counting measures%
\footnote{
The counting measure of a set is the number of elements in the set.}.
%
%
The density \,$ f $\, associated with the measure \,$ \pi $\, consists then in the (discrete) collection of numbers \,$ f_i = f(x_i) $\, such that, for every set \,$ E \subseteq X $\,,
\,$ \pi[E] \, = \, \int_E f \, d\mu \, = \, \sum_{x_i\in E} f(x_i) $\,,
while the density \,$ g $\, (that we may denote \,$ g = \varphi[\,f\,] $\,) associated to the measure \,$ \tau = \varphi[\pi] $\, consists in the (discrete) collection of numbers \,$ g_\alpha = g(y_\alpha) $\, such that
for every set \,$ F \subseteq Y $\, 
\,$ \tau[F] \, = \, \int_F g \, d\nu \, = \, \sum_{y_\alpha\in F} g(y_\alpha) $\,.
Some easy computations then provide the solution:
\begin{equation}
g(y_\alpha) \ = \ \sum_{x_i \in \varphi^\mone[\{ y_\alpha \}]} f(x_i) \qquad \text{for every} \ \ y_\alpha \in Y \quad .
\end{equation}

%
%
\subsection{Propagation of uncertainties in physical measurements}
\label{sec: Propagation of uncertainties in physical measurements}

Physical quantities are often defined in terms of other physical quantities. For instance, the electric resistance \,$ R $\, of a wire is defined as the ratio of the voltage \,$V $\, applied to the wire and the current intensity \,$ I $\, flowing in the wire. Then, a typical measure of \,$ R $\, involves, in fact, the measure of the two quantities \,$ V $\, and \,$ I $\ and the computation of the ratio \,$ R = V/I $\,. 

So, more generally, when one wants to perform a physical measurement of the value of some physical quantity, say \,$ y $\,, most of the time, one resorts to measuring in fact some other quantities, say \,$ \{x^1,x^2,\dots,x^p\} $\,, and then one computes the value of \,$ y $\, via its definition
\begin{equation}
y \ = \ \varphi(x^1,x^2,\dots,x^p) \quad .
\label{eq: mapf1_8361}
\end{equation}
One very basic problem in metrology is that of ``propagating'' the uncertainties appearing in the measurements of the quantities \,$ x^i $\, into the uncertainty on the quantity \,$ y $\,.
Good metrology practice corresponds (ISO, 1993) to representing the uncertainties on a measurement by a probability density (as opposed to simple ``uncertainty bar'').
Therefore, one faces the following problem:

\medskip

\emph{Question:} 
One has some probability measure \,$ \pi $\, defined on the quantities \,$ \{x^1,x^2,\dots,x^p\} $\,, and one \emph{defines} the quantity \,$ y $\, via the mapping in equation~\eqref{eq: mapf1_8361}. What probability measure \,$ \tau $\,  does this imply on the quantity~\,$ y $\,?

\medskip

\emph{Short answer:} The probability measure \,$ \tau $\, is the image of the probability measure \,$ \pi $\,, i.e., according to our general definition of image of a measure: \,$ \tau = \varphi[\pi] $\,.

\medskip

But let us state the problem using a more general terminology.

\medskip

\emph{Preliminaires:}
Some of the quantities \,$ x = \{ x^1,x^2,\dots,x^p \} $\, may be discrete, while others may be real quantities, 
each taking values inside some interval (open or closed). Let \,$ X $\, be the set (part discrete, part continuous) 
whose elements correspond to all the possible values of the quantities \,$ x $\,. Introducing an appropriate $\sigma$-field of subsets of \,$ X $\, is, generally, quite easy%
\footnote{
A Cartesian product of some Borel fields ---for the real variables--- times the collections of all the possible subsets ---for the discrete variables---.},
%
%
so one immediately faces a measurable space \,$ (X,\calE) $\,. We can consider, for more generality, that the quantity \,$ y $\, also is ``multidimensional'': \,$ y = \{y^1,y^2,\dots,y^q\} $\,. A measurable space \,$ (Y,\calF) $\, is introduced as above. Unless that mapping \,$ x \mapsto y = \varphi(x) $\, is pathological%
\footnote{
Physicists need to try hard before being able to introduce mappings that are not measurable with respect to the obvious topologies.},
%
%
it will be measurable (with respect the two $\sigma$-fields \,$ \calE $\, and \,$ \calF $\,).
The (uncertain) result of the measurement of the quantities \,$ x $\, is represented by a $\sigma$-finite%
\footnote{
Physicists will typically represent their measurement uncertainties by introducing probability densities and
discrete probabilities, to be interpreted as the Radon-Nikodym derivatives of the measure \,$ \pi $\,, so 
\,$ \pi $\, will be $\sigma$-finite by construction.}
%
%
measure \,$ \pi $\, on \,$ (X,\calE) $\,.

\medskip

\emph{Question:}
How do the uncertainties encapsulated by the measure \,$ \pi $\, ``propagate'' into uncertainties on the space \,$ (Y,\calF) $\,,
i.e., which is the measure, say \,$ \tau $\,, implied on \,$ (Y,\calF) $\, by the measure \,$ \pi $\, and the mapping \,$ \varphi $\,?

\medskip

\emph{Answer:}
The notion of ``propagation of uncertainties'' can be made precise by imposing that the probability \,$ \tau[B] $\, 
of any subset \,$ B \in \calF $\, must equal the probability of the pre-image (or reciprocal image) of the subset:
\begin{equation}
\tau[B] \ = \ \pi[\,\varphi^\mone [B]\,] \qquad \text{for any} \ B \in \calB \quad .
\end{equation}
But this is exactly our definition of image of a measure, so the answer is
\begin{equation}
\tau \ = \ \varphi[\pi] \quad .
\end{equation}

\medskip

\medskip

\emph{Example:}
The measurement of an electric resistance \,$ R $\, involves the measurement of the two quantities \,$ V $\, and \,$ I $\, and the use of the definition \,$ R = V/I $\,. If the result of the measurement of \,$ V $\, and \,$ I $\, (and the associated uncertainties) is that represented by the (lognormal) probability density
\begin{equation}
f(V,I) \ = \ \frac{1}{2 \, \pi \, \sigma_V \, \sigma_I} \frac{1}{V \, I} \exp\Big(-\frac{\log^2(V/V_0)}{2 \, \sigma_V^2} - \frac{\log^2(I/I_0)}{2 \, \sigma_I^2} \Big) \quad ,
\label{eq: braz_947723}
\end{equation}
then, the notion of image of a measure produces, for the electric resistance \,$ R $\,, the (lognormal) probability density
\begin{equation}
g(R) \ = \ \frac{1}{\sqrt{2 \, \pi} \, \sigma_R} \frac{1}{R} \exp\Big(-\frac{\log^2(R/R_0)}{2 \, \sigma_R^2} \Big) \quad ,
\label{eq: braz_947724}
\end{equation}
where
\begin{equation}
R_0 \ = \ V_0/I_0 \qquad \text{and} \qquad \sigma_R \ = \ \sqrt{\sigma_V^2 + \sigma_I^2} \quad .
\end{equation}
Let us see some details of that.
The space \,$ X $\, is \,$ (0,\infty) \times (0,\infty) \subset \Re^2 $\,, that we endow with two coordinates \,$ \{V,I\} $\, (having the physical interpretation of an electric voltage and an electric intensity). The space \,$ Y $\, is \,$ (0,\infty) \subset \Re$\, that we endow with a coordinate \,$ R $\, (having the physical interpretation of an electric resistance). The mapping \,$ \varphi $\, is (definition of electric resistance) \,$ R = V/I $\,. The usual Borel $\sigma$-fields of \,$ X $\, and of \,$ Y $\, (say \,$ \calB_2 $\, and \,$ \calB_1$\,) are introduced, and the usual Lebesgue measures are considered as base measures.
To arrive at the density \,$ g(R) $\, one can here introduce the ``slack'' variable \,$ P = V \, I $\,, this allowing to consider the ``change of variables'' \,$ \{V,I\} \mapsto \{ R,P \} $\,. One then easily evaluates the density \,$ g_0(R,P) $\, (using the Jacobian of the transformation), and, from it,
\,$ g(R) = \int_0^\infty g_0(R,P) \, dP $\,. It can be shown that the final result for \,$ g(R) $\, is independent on the particular choice of slack variable.

%
%
\subsection{Interpretation of observations (1: using sets)}
\label{sec: Interpretation of observations (1: using sets)}

In the physical sciences, some problems of interpretation of observations can be idealized as follows.
There are two sets \,$ X $\, and \,$ Y $\,, a mapping \,$ \varphi $\, from \,$ X $\, into \,$ Y $\,, and

\medskip

\noindent
\emph{(i)}
we are interested in identifying a particular element \,$ x \in X $\,,
and we have the ``a~priori information'' that it belongs to a subset
\,$ X_\text{prior} \subseteq X $\,:
\begin{equation}
x \ \in \ X_\text{prior} \quad ,
\end{equation}
\emph{(ii)}
we have ``observed'' that some element \,$ y \in Y $\, belongs 
to a subset \,$ Y_\text{obs} \subseteq Y $\,:
\begin{equation}
y \ \in \ Y_\text{obs} \quad ,
\end{equation}
and \emph{(iii)} we know that \,$ y $\, is related to \,$ x $\, via the mapping
\,$ \varphi $\,:
\begin{equation}
y \ = \ \varphi(x) \quad .
\end{equation} 

\begin{figure}[htbp]
   \centering
   \includegraphics[width = 120 mm]{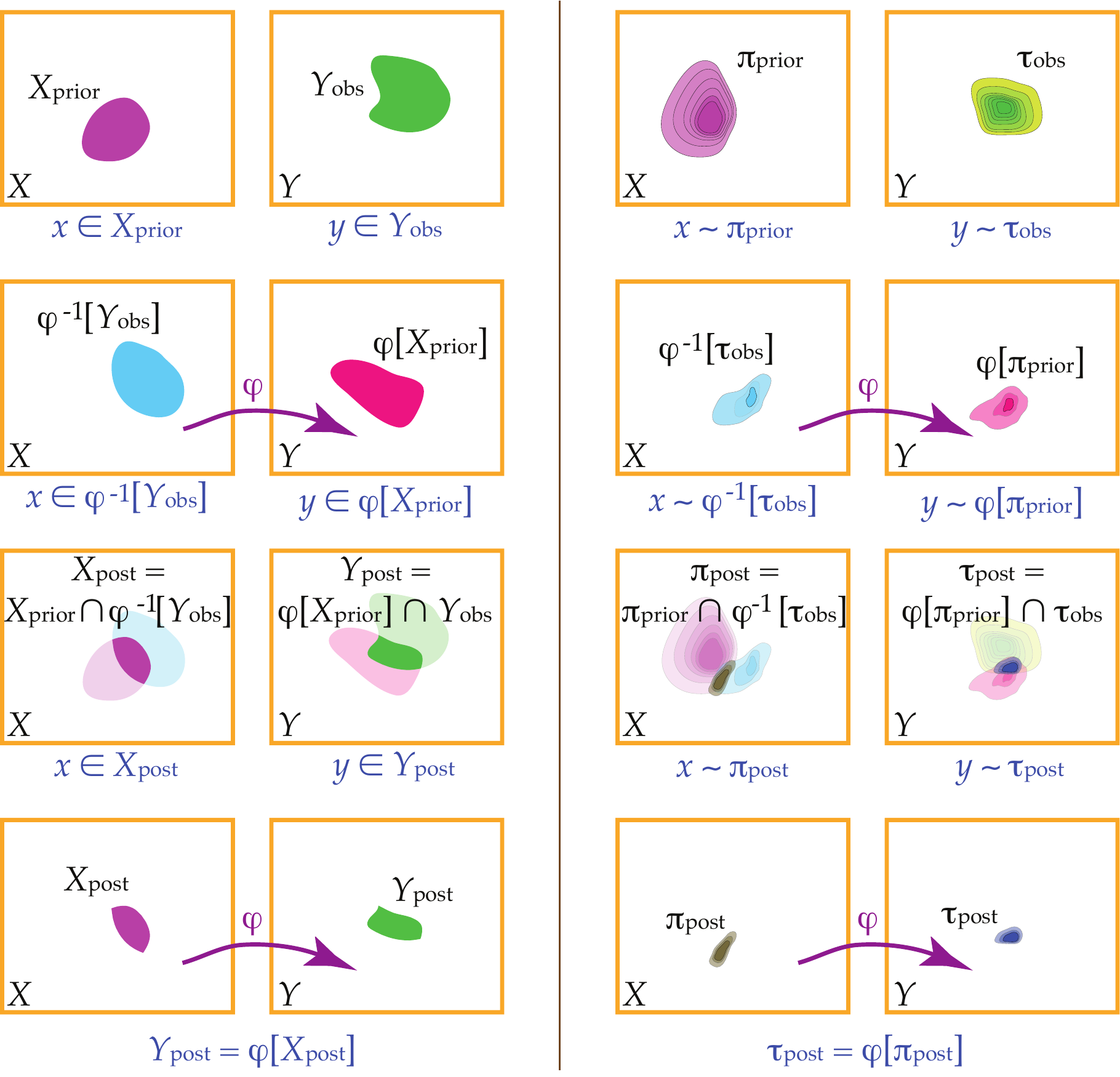} 
   \caption{At the left, an inference problem that can be solved using only set theory (see text). 
   At the right, a similar problem, but this one concerning probability measures (see text). 
   This example at the right corresponds to many of the so-called
   \emph{inverse problems} in the experimental sciences.}
   \label{fig: Mappings}
\end{figure}

These three pieces of information, when put together (see the left of figure~\ref{fig: Mappings}), allow one to infer (using standard set theory reasoning):

\medskip

\noindent
\emph{(i)}
that the element \,$ x $\, belongs, in fact, to a set \,$ X_\text{post} $\, that is smaller or equal to the original set \,$ X_\text{prior} $\,,
\begin{equation}
x \ \in \ X_\text{post} \ = \ X_\text{prior} \cap \varphi^\mone [Y_\text{obs}] \ \subseteq \ X_\text{prior} 
\quad ,
\end{equation}
\noindent
\emph{(ii)}
while the element \,$ y $\, belongs, in fact, to a set \,$ Y_\text{post} $\, that is smaller or equal to the original set \,$ Y_\text{obs} $\,,
\begin{equation}
y \ \in \ Y_\text{post} \ = \ \varphi[X_\text{prior}] \cap Y_\text{obs} \ \subseteq \ Y_\text{obs}
\quad .
\end{equation}

\medskip

These two results are obvious.
Perhaps less obvious is the relation
\begin{equation}
\boxed{ \qquad
Y_\text{post} \ = \ \varphi[X_\text{post}] \quad , \quad 
} 
\end{equation}
that follows directly from the universal set property 
\,$ \varphi[ \, A \cap \varphi^\mone [ \, B \, ] \, ] \, 
= \, \varphi[A] \cap B $\, 
(equation~\eqref{eq: forisaad_332}).

Remark that we are inside the paradigm typical of a ``problem of 
assimilation of observations'' ---sometimes called 
``inverse modeling problem''---: the mapping \,$ x \mapsto y = \varphi(x) $\, can be seen as the typical mapping between the ``model parameter space'' and the ``observable parameter space''. 
In what concerns the element \,$ x \in X $\, 
we pass from the ``a priori information''
\,$ x \in X_\text{prior} $\, to the 
``a posteriori information'' \,$ x \in X_\text{post} \subseteq X_\text{prior} $\,. Similarly, in what concerns
the element \,$ y \in Y $\, we pass from the ``initial observation'' 
\,$ y \in Y_\text{obs} $\, to the ``refined observation'' 
\,$ y \in Y_\text{post} \subseteq Y_\text{obs} $\,.

In the next section, the same problem is reformulated, but using probability measures instead of sets.

%
%

%
%

%
%
\subsection{Interpretation of observations (2: using measures)}
\label{sec: Interpretation of observations (2: using measures)}

The problem of interpretation of observation ---sometimes called the ``inverse problem''-- appears as follows.
A physical system (e.g., a molecule, an ocean, a~planet's atmosphere, a galaxy) is under investigation. For the purposes of the investigation, the system is described using a collection of \,$ p $\, physical quantities \,$ x = \{x^1,x^2,\dots,x^p\} $\,; some of them taking only discrete values (for instance, \emph{black} or \emph{white}) and some others taking continuous values (for instance a temperature can take any positive real value). A set \,$ X $\, is introduced, the elements (or ``points'') of which corresponding to the quantities \,$ x $\, taking all their possible values. In the jargon of inverse problem theory, the set \,$ X $\, is called the \emph{model space}%
\footnote{
In fact, the set \,$ X $\, is something more abstract: any (invertible) change of variables \,$ x \rightleftharpoons x' $\, is to be seen as a ``change of coordinates'' inside \,$ X $\,, not as the definition of a new set \,$ X' $\,. For a discussion of this kind of intrinsic view on physical quantities, see Tarantola (2006).}.
%
%
In order to gain information on the \,$ p $\, physical quantities \,$ x $\,, a~set of \,$ q $\, physical quantities \,$ y = \{y^1,y^2,\dots,y^q\} $\, ---perhaps only quite indirectly related to the quantities \,$ x $\,--- is measured. As above, when considering all possible values for the \,$ q $\, quantities \,$ y $\, one is faced with a set \,$ Y $\,, the ``observable parameter space''. One then (usually implicitly) considers two collections of subsets \,$ \calE $\, and \,$ \calF $\,, of sets of \,$ X $\, and of \,$ Y $\, respectively, that are $\sigma$-algebras, so one has two measurable spaces \,$ (X,\calE) $\, and \,$ (Y,\calF) $\,. On each of these spaces, a base ($\sigma$-finite) measure has to be considered in order to have two $\sigma$-finite measure spaces, \,$ (X,\calE,\mu) $\,, and \,$ (Y,\calF,\nu) $\,. To a mathematician, the existence of the two base measures \,$ \mu $\, and \,$ \nu $\, may seem a minor hypothesis. A physicist may have to work hard to find them, as they must represent the volume measure of each space, and, as such, they must have the necessary invariances%
\footnote{
For instance, one of the quantities may be the period \,$ T $\, (of, say, a star). How to measure the volume (in fact, the length) of an interval \,$ (T_1,T_2) $\,, say \,$ \mu[(T_1,T_2)] $\,? Taking \,$  \mu[(T_1,T_2)] = |T_2-T_1| $\, would not be consistent with, when working with the frequency \,$ \omega=2\pi/T $\,, measuring the volume as \,$ |\omega_2-\omega_1| $\,, because \,$ |T_2-T_1| \neq |\omega_2-\omega_1| $\,. In fact, the right volume measure is \,$ \mu[(T_1,T_2)] = |\log(T_2/T_1)| $\,, because \,$ |\log(T_2/T_1)| = |\log(\omega_2/\omega_1)| $\,. See Tarantola (2005) for an elementary discussion of this problem, or Tarantola (2006) for a more advanced discussion.}.
%
%
For this reason, let us here call the two base measures \,$ \mu $\, and \,$ \nu $\, the respective \emph{volume measures}. They matter, because the reciprocal image of a probability measure on \,$ (Y,\calF) $\, and the intersection of measures on \,$ (X,\calE) $\, and on \,$ (Y,\calF) $\, depend on them.

The final structure element is that a physical theory is assumed to exist, that is able ---given any possible value of the model parameters \,$ x $\,--- to predict (in a Popperian sense) the observations \,$ y $\,. This prediction consists in a mapping \,$ \varphi : X \mapsto Y $\,, that must be assumed to be measurable (what, for a physicist, just means that \,$ \varphi $\, is assumed to be not ``pathological''). Of course, the mapping \,$ \varphi $\, is not assumed to be invertible (and it may be ``nonlinear'').

In a typical inverse problem one cares in introducing any available a priori information on \,$ x $\, 
(that means information available before the measurements on \,$ y $\, are carried out) as a probability measure, say \,$ \pi_\text{prior} $\,, on \,$ (X,\calE) $\, (it must be ``ordinary'', i.e., $\sigma$-finite, but it does not need to be absolutely continuous w.r.t.\ the volume measure \,$ \mu $\,). When the measurements of the quantities \,$ y $\, are carried out, the result is represented%
\footnote{
The representation of the (possibly uncertain) result of an observation as a probability measure is in compliance with ISO's (1993) recommendations.}
%
%
as a probability measure, say \,$ \tau_\text{obs} $\,, on \,$ (Y,\calF) $\,, that must be absolutely continuous w.r.t.\ the volume measure \,$ \nu $\, (i.e., it has to be representable by a density).
So, one has the following three elements (see the right of figure~\ref{fig: Mappings}):

\medskip

\noindent
\emph{(i)}
a priori information on the model parameters, i.e., a probability measure on \,$ (X,\calE) $\,
\begin{equation}
\pi_\text{prior} \quad ;
\end{equation}
\emph{(ii)}
results of the measurements, i.e., a probability measure on \,$ (Y,\calF) $\,
\begin{equation}
\tau_\text{obs} \quad ;
\end{equation}
and \emph{(iii)}
the modeling mapping
\begin{equation}
\varphi: X \mapsto Y \quad .
\end{equation}

\medskip

It is clear that the existence of a mapping \,$ \varphi $\, is going to transform the measure \,$ \pi_\text{prior} $\, into some other measure, say \,$ \pi_\text{post} $\,, and the measure \,$ \tau_\text{obs} $\, into some measure, say \,$ \tau_\text{post} $\,, much as it happened in section~\ref{sec: Interpretation of observations (1: using sets)}, where the prior sets were transformed into posterior sets. The problem is that here we are facing natural objects (measurements and physical laws), that are not easily amenable to axiomatization. Usual presentations of the inverse problem unconvincingly use intuitive interpretations of the notion of conditional probability and of, perhaps, Bayes' theorem. I prefer here to frankly state that I formulate the problem using only the analogy between the present problem and the set-theoretical problem in section~\ref{sec: Interpretation of observations (1: using sets)} (although closer analogies can be elaborated%
\footnote{
Assume that a random \,$ x $\, and a random \,$ y $\, are created according respectively to \,$ \pi_\text{prior} $\, and \,$ \tau_\text{obs} $\,, and that the pair \,$ \{x,y\} $\, is accepted only if \,$ y = \varphi(x) $\, 
(much as we did in section~\ref{sec: Intersection of probability measures}). It is easy to prove (see the argument in section~\ref{sec: Intersection of probability measures}) that when a pair \,$ \{x,y\} $\, is accepted, \,$ x $\,~is a sample point of the measure \,$ \pi_\text{post} = \pi_\text{prior} \cap \varphi^\mone[\tau_\text{obs}] $\, and \,$ y $\, is a sample point of the measure \,$ \tau_\text{post} = \varphi[\pi_\text{prior}] \cap \tau_\text{obs} $\,. These are exactly expressions~\eqref{eq: merdefaitbeau_772} and~\eqref{eq: merdefaitbeau_773}.}).
%
%
The results there (that concerned intersection of sets, and images and reciprocal images of sets)
were unquestionable. As far as the present theory, defining the intersection of measures, and images and reciprocal images of measures is an acceptable generalization of (a part of) set theory (and the compatibility property suggests that it is), we can match this problem to the one in section~\ref{sec: Interpretation of observations (1: using sets)}
(see also the parallel suggested in figure~\ref{fig: Mappings}). In any case, the formulas we are going to find for the inverse problem are basically identical (although, perhaps, a little more general) than those proposed in the usual literature%
\footnote{
For a probabilistic formulations of the inverse problem, see Tarantola and Valette (1982), Menke (1989), Mosegaard and Tarantola (1995), Aster et al.\ (2005), or Tarantola (2005). For an alternative, statistical decision theory, see Evans and Stark (2002).}.
%
%

Then (see the illustration at the right of figure~\ref{fig: Mappings}):

\medskip

\noindent
\emph{(i)}
on the model parameter space \,$ (X,\calE) $\,, one passes from the prior probability measure \,$ \pi_\text{prior} $\, to the posterior probability measure
\begin{equation}
\pi_\text{post} \ = \ \pi_\text{prior} \cap \varphi^\mone[\tau_\text{obs}] \quad ;
\label{eq: merdefaitbeau_772}
\end{equation}
\emph{(ii)}
on the observable parameter space \,$ (Y,\calF ) $\,, one passes from the initial probability measure \,$ \tau_\text{obs} $\, (representing the result of the measurements) to the probability measure
\begin{equation}
\tau_\text{post} \ = \ \varphi[\pi_\text{prior}] \cap \tau_\text{obs} 
\label{eq: merdefaitbeau_773}
\end{equation}
representing a refined estimation of the values of the observable parameters \,$ y $\,.
Finally, \emph{(iii)} the compatibility property (equation~\eqref{eq: forisaad}) nicely states that
\begin{equation}
\tau_\text{post} \ = \ \varphi[\pi_\text{post}] \quad .
\label{eq: merdefaitbeau_774}
\end{equation}

Let us evaluate the posterior probability \,$ \pi_\text{post}[E] $\, of some set \,$ E \in \calE $\,. 
From expression~\eqref{eq: merdefaitbeau_772} it follows (using, first, the definition of the intersection of measures in equation~\eqref{eq: intersec_3822}, then, the definition of the reciprocal image of a measure in equation~\eqref{eq: erecp_46})
\begin{equation}
\pi_\text{post}[E] \ = \ \frac{1}{n} \int_E \frac{d\tau}{d\nu}(\varphi(x)) \ \, d\pi_\text{prior} \quad ,
\label{eq: posterior_9924}
\end{equation}
where the constant \,$ n = \int_X \frac{d\tau}{d\nu}(\varphi(x)) \, d\pi_\text{prior} $\, must be different from zero (in order for \,$ \pi_\text{post} $\, to be a probability measure). Note that the density \,$ d\tau/d\nu $\, exists because \,$ \tau $\, was assumed to be absolutely continuous w.r.t.\ the volume measure \,$ \nu $\,.

To evaluate the posterior probability \,$ \tau_\text{post}[F] $\, of a set \,$ F \in \calF $\,, one could try to start with expression~\eqref{eq: merdefaitbeau_773}, but this possibility is not the most practical. One can rather use the compatibility relation~\eqref{eq: merdefaitbeau_774}, as then (because of the relation~\eqref{eq: erecp_46} defining the image of a measure),
\begin{equation}
\tau_\text{post}[F] \, = \, \pi_\text{post}[\varphi^\mone[F]] \quad .
\label{eq: posterior_9924543}
\end{equation}

In real-life problems, the finite probabilities \,$ \pi_\text{post}[E] $\, and \,$ \tau_\text{post}[F] $\, (i.e., the sums in equations~\eqref{eq: posterior_9924}--\eqref{eq: posterior_9924543} can (approximately) be evaluated using Monte Carlo methods%
\footnote{
From expression~\eqref{eq: posterior_9924} is follows that if \,$ x_1,x_2,\dots $\, is a collection of (independent) random sample elements of the prior probability measure \,$ \pi_\text{prior} $\,, and if, for every element, a~random decision is taken to conserve or discard it with the probability of being conserved equal to \,$ k \, \frac{d\tau}{d\nu}(\varphi(x_i)) $\, (where the positive constant \,$ k $\, is arbitrary, excepted that it must ensure that the maximum attained value is $\leq 1$\,), then, the collection \,$ x_1', x_2', \dots, $\, of conserved elements is a sample of \,$ \pi_\text{post} $\, (Mosegaard and Tarantola [1995]). And (as it follows from the definition of image of a measure) the collection \,$ \varphi(x'_1) , \varphi(x'_2) , \dots $\, is a sample of \,$ \tau_\text{post} $\,.}.

As a final remark, should the prior probability measure \,$ \pi_\text{prior} $\, be absolutely continuous w.r.t.\ the volume measure \,$ \mu $\, (i.e., should the density \,$ d\pi_\text{prior}/d\mu $\, exist), the posterior probability measure \,$ \pi_\text{post} $\, would also have a density, whose explicit expression would follow immediately from equation~\eqref{eq: posterior_9924}:
\begin{equation}
\frac{d\pi_\text{post}}{d\mu} \ = \ \frac{1}{n} \ \frac{d\pi_\text{prior}}{d\mu} \ \frac{d\tau}{d\nu} \circ \varphi \quad ,
\end{equation}
i.e., explicitly,
\,$ \frac{d\pi_\text{post}}{d\mu}(x) \ = \ \frac{1}{n} \, \frac{d\pi_\text{prior}}{d\mu}(x) \, \frac{d\tau}{d\nu}(\varphi(x)) $\,. 
In the jargon of inverse theory, \,$ x \mapsto \frac{d\tau}{d\nu}(\varphi(x)) $\, is called the ``likelihood function''.

%
%

\section{References}

\bibref
Ambrosio, L., N.\ Gigli, and G.\ Savar\'e, 2005. Gradient flows: in metric spaces and in the space of probability measures, {\sl Birk\"auser Basel}.

\bibref
Aster, R.C., C.H.\ Clifford, and B.\ Borchers, 2005. Parameter estimation and inverse problems, {\sl Academic Press}.

\bibref
Evans, S.N.\ and P.B.\ Stark, 2002. Inverse problems as statistics, {\sl Inverse problems}, {\bf 18}, R55--97.

\bibref
Halmos, P.R., 1950. Measure theory, {\sl Springer-Verlag}.

\bibref
International Organization for Standardization (ISO), 1993. Guide to the expression of uncertainty in measurement, {\sl ISO}.

\bibref
Kolmogorov, A.N., 1950. Foundations of the theory of probability, {\sl Chelsea publishing company}.

\bibref
Kuo, H.-H., 2002, White noise theory, {\sl in:} Handbook of stochastic analysis and applications (editors: Kannan and Lakshmikantham), {\sl CRC Press}. 

\bibref
Menke, W., 1989. Geophysical data analysis: discrete inverse theory, {\sl Academic Press}.

\bibref
Mosegaard, K., and A\ Tarantola, 1995. Monte Carlo sampling of solutions to inverse problems, 
{\sl J.\ Geophys.\ Res.}, Vol.\ 10, no.\ B7, p.\ 12,431--12,447.

\bibref
Rudin, W., 1970. Real and complex analysis, {\sl McGraw-Hill}.

\bibref
Schilling, R.L., 2006. Measures, integrals and martingales, {\sl Cambridge University Press}.

\bibref
Tarantola, A., and Valette, B., 1982. Inverse problems = quest for information. {\sl J. geophys.}, 50, p.\ 159--170.

\bibref
Tarantola, A., 2005. Inverse problem theory and methods for model parameter estimation, {\sl SIAM}.

\bibref
Tarantola, A., 2006. Elements for physics; quantities, qualities, and intrinsic theories, {\sl Springer}.

\bibref
Taylor, S.J., 1966. Introduction to measure theory and integration, {\sl Cambridge University Press}.

\bibref
Weinberg, S., 1972. Gravitation and cosmology, {\sl Wiley}.

\bibref
Yokoi, Y., 1990. Positive generalized white noise functionals, {\sl Hiroshima Math.~J.}, vol.\ 20, no.\ 1, pages 137--157. 

%
%

\section{Acknowledgements}

My first thanks go to Olivier Pironneau, who has constantly guided my hand in the right direction. Without his help, the compatibility conjecture may never have become a theorem. Also, a lot of thanks to Bartolom\'e Coll, for systematically guiding my hand in the wrong direction. His intelligence and patience are legendary, as is his obstinacy: this forced me to really understand (I hope) what I was trying to do. Klaus Mosegaard is the accomplice of many of my adventures. It was a question of him ({\sl Albert, are you certain that there are no simpler formulas for stating the general solution of an inverse problem?}) that prompted me to abandon the use of conditional probabilities, which resulted in the making of the present theory. Philip Stark has a view of probabilistic reasoning that is almost opposite to mine. But I have learned from him how to use probability theory rigorously. His e-mail with a few pages of Halmos' book was crucial in solving my last difficulties with the theory. Finally, many thanks to Bob Nowack for a critical reading of the manuscript.

\end{document}